\documentclass[11pt]{amsart}
\usepackage{amssymb,amsmath,epsfig,mathrsfs, enumerate}
\usepackage{graphicx}
\usepackage[normalem]{ulem}
\usepackage{fancyhdr}
\usepackage[margin=3cm]{geometry}
\usepackage{hyperref}
\hypersetup{
 colorlinks   = true,
 urlcolor     = blue,
 linkcolor    = blue,
 citecolor   = red ,
 bookmarksopen=true
}

\theoremstyle{plain}
\newtheorem{thrm}{Theorem}[section]
\newtheorem{lemma}[thrm]{Lemma}

\newtheorem{rmrk}[thrm]{Remark}
\newtheorem{dfn}[thrm]{Definition}

\begin{document}

\newcommand{\SL}{\mathcal L^{1,p}( D)}
\newcommand{\Lp}{L^p( Dega)}
\newcommand{\CO}{C^\infty_0( \Omega)}
\newcommand{\Rn}{\mathbb R^n}
\newcommand{\Rm}{\mathbb R^m}
\newcommand{\R}{\mathbb R}
\newcommand{\Om}{\Omega}
\newcommand{\Hn}{\mathbb H^n}
\newcommand{\aB}{\alpha B}
\newcommand{\eps}{\ve}
\newcommand{\BVX}{BV_X(\Omega)}
\newcommand{\p}{\partial}
\newcommand{\IO}{\int_\Omega}
\newcommand{\bG}{\boldsymbol{G}}
\newcommand{\bg}{\mathfrak g}
\newcommand{\bz}{\mathfrak z}
\newcommand{\bv}{\mathfrak v}
\newcommand{\Bux}{\mbox{Box}}
\newcommand{\e}{\ve}
\newcommand{\X}{\mathcal X}
\newcommand{\Y}{\mathcal Y}
\newcommand{\W}{\mathcal W}
\newcommand{\la}{\lambda}
\newcommand{\vf}{\varphi}
\newcommand{\rhh}{|\nabla_H \rho|}
\newcommand{\Ba}{\mathcal{B}_\beta}
\newcommand{\Za}{Z_\beta}
\newcommand{\ra}{\rho_\beta}
\newcommand{\na}{\nabla_\beta}
\newcommand{\vt}{\vartheta}

\numberwithin{equation}{section}

\newcommand{\RN} {\mathbb{R}^N}
\newcommand{\Sob}{S^{1,p}(\Omega)}
\newcommand{\Dxk}{\frac{\partial}{\partial x_k}}
\newcommand{\Co}{C^\infty_0(\Omega)}
\newcommand{\Je}{J_\ve}
\newcommand{\beq}{\begin{equation}}
\newcommand{\bea}[1]{\begin{array}{#1} }
\newcommand{\eeq}{ \end{equation}}
\newcommand{\ea}{ \end{array}}
\newcommand{\eh}{\ve h}
\newcommand{\Dxi}{\frac{\partial}{\partial x_{i}}}
\newcommand{\Dyi}{\frac{\partial}{\partial y_{i}}}
\newcommand{\Dt}{\frac{\partial}{\partial t}}
\newcommand{\aBa}{(\alpha+1)B}
\newcommand{\GF}{\psi^{1+\frac{1}{2\alpha}}}
\newcommand{\GS}{\psi^{\frac12}}
\newcommand{\HFF}{\frac{\psi}{\rho}}
\newcommand{\HSS}{\frac{\psi}{\rho}}
\newcommand{\HFS}{\rho\psi^{\frac12-\frac{1}{2\alpha}}}
\newcommand{\HSF}{\frac{\psi^{\frac32+\frac{1}{2\alpha}}}{\rho}}
\newcommand{\AF}{\rho}
\newcommand{\AR}{\rho{\psi}^{\frac{1}{2}+\frac{1}{2\alpha}}}
\newcommand{\PF}{\alpha\frac{\psi}{|x|}}
\newcommand{\PS}{\alpha\frac{\psi}{\rho}}
\newcommand{\ds}{\displaystyle}
\newcommand{\Zt}{{\mathcal Z}^{t}}
\newcommand{\XPSI}{2\alpha\psi \begin{pmatrix} \frac{x}{|x|^2}\\ 0 \end{pmatrix} - 2\alpha\frac{{\psi}^2}{\rho^2}\begin{pmatrix} x \\ (\alpha +1)|x|^{-\alpha}y \end{pmatrix}}
\newcommand{\Z}{ \begin{pmatrix} x \\ (\alpha + 1)|x|^{-\alpha}y \end{pmatrix} }
\newcommand{\ZZ}{ \begin{pmatrix} xx^{t} & (\alpha + 1)|x|^{-\alpha}x y^{t}\\
     (\alpha + 1)|x|^{-\alpha}x^{t} y &   (\alpha + 1)^2  |x|^{-2\alpha}yy^{t}\end{pmatrix}}
\newcommand{\norm}[1]{\lVert#1 \rVert}
\newcommand{\ve}{\varepsilon}

\title[space like unique continuation  etc.]{Space like strong unique continuation for sublinear parabolic equations}

\author{Agnid Banerjee}
\address{Tata Institute of Fundamental Research\\
Centre For Applicable Mathematics \\ Bangalore-560065, India} \email[Agnid Banerjee]{agnidban@gmail.com}
%
%
%
%
%
\author{Ramesh Manna}
\address{Tata Institute of Fundamental Research\\
Centre For Applicable Mathematics \\ Bangalore-560065, India}
\email[Ramesh Manna]{ramesh@tifrbng.res.in}
%


\thanks{Second author supported by  SERB National Postdoctoral fellowship, PDF/2017/0027}


%
%
%

\maketitle

\tableofcontents

\begin{abstract}

In this paper, we establish space like strong unique continuation property  (sucp) for uniformly parabolic sublinear equations under appropriate structural assumptions. Our main result Theorem \ref{main} constitutes the parabolic counterpart  of  the strong unique continuation result recently established in \cite{Ru} for analogous elliptic sublinear equations.  Similar to that in \cite{Ru},  this is accomplished  via a   new  $L^{2}-L^{2}$ type Carleman estimate for a class of sublinear parabolic operators. 
\end{abstract}

\section{Introduction}
\noindent     The primary objective of this paper is to study space like strong unique continuation for backward sublinear second order parabolic operators as in \eqref{sub} below with structural assumptions on the sublinearity as in \eqref{a1}.   To begin with, we note  that an  operator $L$  (local or non-local) is said to possess the  strong unique continuation property  if  any non-trivial solution $u$ to
\[
Lu=0
\]
in a (connected) domain $\Om \subset \mathbb{R}^n$ cannot vanish to infinite order   at any point in $\Om$.   An operator $L$  instead is said to  have the weak unique continuation property (wucp)  if a non-trivial solution to 
$Lu=0$  cannot vanish in an open subset. Likewise,  space-like strong unique continuation property for a parabolic operator $L$ asserts that if a solution $u$ to 
  \[
  Lu=0
  \]
   vanishes to infinite order at some point $(x_0, t_0)$,  then  $u(\cdot, t_0) \equiv 0$. The unique continuation property for second order elliptic and parabolic equations has a long history  and by now has several important ramifications.

\medskip

A prototypical example of an operator $L$ which satisfies the strong unique continuation property is the Laplacian $\Delta$ in which case  sucp is a consequence of the real analyticity of solutions to
\[
\Delta u=0.
\]
This property is however true for more general elliptic equations of the type
\begin{equation}\label{var}
\operatorname{div}(A(x) \nabla u) + b.\nabla u +  Vu=0
\end{equation}
where the principal part $A$ can be allowed to be Lipschitz and where $b, V$ have appropriate  integrability properties.   Based on a visionary work due to Carleman in 1939 ( see \cite{Car}) who established strong unique continuation for
\[
-\Delta + V, \ V \in L^{\infty}
\]
in $\R^2$, Carleman  estimates were developed systematically in the seminal work of \cite{AKS} in 1962 where sucp was established for Lipchitz $A$ and bounded $b,V$.  We note that Lipschitz regularity assumption on $A$ is optimal in view of a deep counterexample due to Plis in  \cite{Pl}. Some of the  other important  works in this direction are due to Chanillo-Sawyer (\cite{CS}), Kenig-Ruiz-Sogge (\cite{KRS}) and Jerison-Kenig (\cite{JK}). Each of these works deal with scaling critical potentials in different function spaces. For instance in \cite{JK}, sucp is established for 
\[
-\Delta + V
\]
where $V \in L^{n/2}$. The result of Jerison and Kenig  was subsequently extended  by Koch and Tataru in \cite{KT0} to equations of the type \eqref{var} with borderline Lipschitz principal part. 

\medskip

An alternate approach  which   is  instead based on the almost monotonicity of a generalized   frequency function introduced by Almgren in \cite{Al} came up in the  works of Garofalo and Lin in 1986( \cite{GL1}, \cite{GL2}).  Using this approach, they were able to obtain new quantitative  information on the zero set of solutions to divergence form elliptic equations and in particular,  their  results encompassed that of  \cite{AKS}.   Also in recent times, their approach found application in the optimal regularity of solutions for a class of free boundary problems known as  Signorini problems (see for instance \cite{ACS}, \cite{CSS}).

  \medskip
  
  The study of weak unique continuation for parabolic equations began with the early work of Mizohata in \cite{Mi}   and Yamabe in \cite{Y} followed by the work of Sogge in \cite{So} where  certain classes of unbounded potentials were treated using appropriate $L^{p}$ Carleman estimates.  The study of  strong backward uniqueness  for parabolic  equations  for time independent coefficients began with the work of  Lin in \cite{L}. Subsequently Poon in \cite{Po} established strong backward uniqueness for global solutions  under   Tychonoff type exponential growth assumption on the solution    by  adapting   to the parabolic setting the  frequency function approach of  Garofalo and  Lin.  This continued with the work of Chen in \cite{Ch} where instead Carleman estimates were employed.     We note that backward uniqueness is in general not true without such global assumptions on the solution. This follows from  a counterexample due to Frank Jones in \cite{F} were it is shown that there exists a non-trivial  unbounded caloric function that is supported in a time strip of the type $\R^n \times (t_1, t_2)$.   Moreover, the counterexample of Frank Jones also shows that in general one cannot expect space-time strong unique continuation property for parabolic equations. Therefore  in this scenario,  the question of space like strong unique continuation is more  relevant for local solutions. 
  
\medskip

 This was taken up by Escauriaza and Fernandez in \cite{EF}  where they established space like strong unique continuation property for backward   parabolic equations  of the type
  \begin{equation}\label{pa}
  \operatorname{div}(A(x,t) \nabla u) + \partial_t u  + b(x,t).\nabla  u + Vu =0
  \end{equation}
  where  $b, V$ are bounded  and the principal part $A$ has regularity assumptions  similar to that in \eqref{coef}  below.
  We also refer to  the  subsequent work of Escauriaza-Fernandez-Vessella    in  \cite{EFV} where certain quantitative results were obtained.     The approach in \cite{EF} and \cite{EFV} are based on a $L^{2}-L^{2}$ type Carleman estimate of the type
   \begin{equation}
   ||t^{-\alpha-1/2} e^{-\frac{|x|^2}{8t}} u||_{L^2} \lesssim || t^{-\alpha} e^{-\frac{|x|^2}{8t}} (\partial_t u + \operatorname{div}(A\nabla u))||_{L^2}
   \end{equation}
 which is obtained using a fairly nontrivial parabolic Rellich type identity as stated in Lemma \ref{lma1} below coupled with a clever integration by parts argument.   The  Carleman estimate in \cite{EF}  is also partly  inspired by the previous work of Poon in \cite{Po}. In fact the work of Poon contributed in  clarifying the correct form of Carleman estimates  that can be expected in the parabolic situation.  The space like sucp in \cite{EF} was later on extended in \cite{KT} to parabolic equations   with  principal part $A(x,t)$ having lower regularity   in the time variable  $t$ and where the lower order terms  $b$ and $V$ are allowed  to belong to some scaling critical function spaces similar to that in the elliptic case as in \cite{JK} and \cite{KT0}. We note that unlike that in  \cite{EF}, the  proof in \cite{KT} is  instead based on  deep $L^{p}$ spectral projection bounds  for the Hermite operator.  Such bounds were independently obtained by Thangavelu in \cite{T} and Kharazdhov in \cite{K} and they  were also essential in the proof of $L^{p}$ Carleman estimates for the heat operator in the  previous works of Escauriaza in \cite{E} and Escauriaza and Vega in \cite{EV} using which the authors showed  backward uniqueness for \eqref{pa} when $A=\mathbb{I}$ and $V \in L^{1}L^{\infty} + L^{\infty}L^{n/2}$.

  \medskip

Now regarding  sublinear equations,  we note that   motivated by the study of nonlinear eigenvalue problems, the analysis of corresponding nodal domains as in \cite{PW} and also  because of certain connection of such equations  to porous media type equations (see for instance \cite{Vaz}),   the study of  unique continuation for  sublinear elliptic equations  was taken up in recent times  by Soave and Weth in \cite{SW}  where they established wucp for equations of the type
 \begin{equation}\label{t100}
 \operatorname{div}(A(x) \nabla u) + f(x, u) + Vu=0
 \end{equation}
 where the sublinear term $f$ satisfies the structural  assumptions similar to that in \eqref{a1} below.  Such equations are modeled on
 \begin{equation}\label{s10}
 -\Delta u = |u|^{p-2} u.
 \end{equation}
 Note that  the study of strong unique continuation for  \eqref{s10} cannot be reduced to that for
 \[
 -\Delta + V
 \]
  because in this case, $V= |u|^{p-2}$ need not be in $L^{n/2}$ near the zero set of $u$  as $p \in (1, 2)$.  In fact  such sublinear equations have their intrinsic  difficulties  and this is also partly visible from the fact  that the sign assumption on the sublinearity $f$ in \eqref{a1} is quite crucial because otherwise unique continuation fails. This later fact follows  from  a counterexample in \cite{SW} where it is shown that unique continuation is not true for
 \begin{equation}\label{s100}
 \Delta u= |u|^{p-2} u, \ p \in (1, 2).
 \end{equation}
    In \cite{SW}, the authors adapted the frequency function approach  of Garofalo and Lin and also that of  Garofalo and Smit Vega Garcia as in \cite{GG}. The  question of strong unique continuation for such  sublinear equations  was then later  addressed by Ruland in \cite{Ru} via new Carleman estimates for  the corresponding sublinear elliptic operators.  See  also the  recent interesting work of Soave and Terracini in \cite{ST} where the authors study the following two phase membrane   problem
    \begin{equation}\label{t10}
    -\Delta u= \lambda_{+} (u^{+})^{q-1} - \lambda_{-} (u^{-})^{q-1}, \ \text{where $\lambda_{+}, \lambda_{-} >0$, $q \in [1, 2)$}
    \end{equation}
    and establish a strong unique continuation property as well as  a regularity result for the nodal domains of solutions to such equations.  The key object in their analysis is a monotonicity formula for a $2$-parameter family of Weiss type functionals introduced by Weiss in \cite{We}.   The reader should however  note that although   \eqref{t10}  is  a more general equation  than \eqref{s100}, but  it doesn't encompass the class of equations as  in \eqref{t100}. Therefore the  unique continuation results  in \cite{Ru} and \cite{SW} are not covered  by the results in \cite{ST} and also vice-versa.

    We would also like to refer to a recent work by  two of us with Garofalo  as in \cite{BGR} where  the result of Ruland has been extended to sublinear equations associated to degenerate elliptic Baouendi-Grushin operators $\mathbb{B}_{\gamma}$  defined by
    \[
  \mathbb{B}_{\gamma}= \Delta_z + |z|^{2\gamma} \Delta_t, \ (z,t) \in \R^m \times R^n.
  \]
 The method in \cite{BGR} also slightly simplifies the proof of Ruland  when the principal part is $\Delta$  and moreover our proof of the sublinear Carleman estimate as stated  in \eqref{Car1} below is also inspired in parts by the ideas in \cite{BGR}.  We also note that the recent work \cite{ST}   addresses  the related nodal domain estimates  for solutions to such sublinear equations. 
  
   Therefore given the  recent developments in the sublinear unique continuation theory in the elliptic case as in \cite{SW} and \cite{Ru}, in this paper, we study analogous strong unique continuation for backward  parabolic sublinear equations of the type 
 \begin{equation}\label{sub}
 \sum_{i,j=1}^n \partial_i(a_{ij}(x,t)\partial_ju) +\partial_tu+f(X, u) + Vu= 0
 \end{equation}
 
 where $V \in L^{\infty}$ and  $f$ and  its primitive $F$ satisfies  the following structural  assumptions similar to  that in \cite{Ru} and \cite{SW} for some $\kappa, K > 0$: 
 \begin{equation}\label{a1}
 \begin{cases}
 f((x,t), 0) =0
 \\
 F((x,t), s) = \int_{0}^{s} f(X, s) ds
 \\
  0 < s f((x,t), s) \leq q F((x,t), s), \text{for some $q \in (1, 2)$ and  $s \in (-1, 1) \setminus \{0\}$}
  \\
  |\nabla_{(x,t)} f| \leq K |f|,\ |\nabla_{(x,t)} F| \leq K F
  \\
  f((x,t), s) \leq \kappa s^{p-1}\ \text{for some $p \in (1,2)$};
  \end{cases}
  \end{equation}
  We note that the first  and the last condition  in \eqref{a1} implies that  for some constant $c_0, c_1$, we have that
  \begin{equation}\label{a0}
 c_1 s^{p} \geq  F(\cdot, s) \geq c_0 s^{q}, \text{ for $s \in (-1,1)$}.
  \end{equation}
  A  prototypical $f$ satisfying \eqref{a1}  is given by 
  \[
  f((x,t), u) = \sum_{i=1}^l c_i(x,t) |u|^{q_i-2} u,\] where for each $i$,  $q_i \in (1,2)$,  $0<k_0<c_i< k_1$ and $|\nabla c_i| < K$ 
 for some $k_0, k_1\ \text{and}\ K$.  In this case, we can take  $q= \text{max}\{q_i\}$ and $p= \text{min} \{q_i \}$.  On the principal part $A=(a_{ij})$,  similar to \cite{EF},  we assume that there exists   $C_0 , \la > 0$ and $0< \beta \leq 1 $  such that for all $x$ and $y$ in $\R^n$ and $0\leq t, s< \infty$, we have

\begin{equation}\label{coef}
\begin{cases}
\lambda |\xi|^2 \leq\sum a_{ij} ((x,t), u)  \xi_i \xi_j \leq \lambda^{-1} |\xi|^2,\ \text{for some $\lambda>0$};
\\
|A(x,t) - A(y,s)| \leq C_0(|x-y|^2+|t-s|)^{\beta/2}\ \text{and also}\ |\nabla_x a_{ij}(x,t)| \leq C_0 |x|^{\beta-1}, \ |\partial_t a_{ij}(x,t) | \leq C_0 t^{\beta/2 -1}.
\end{cases}
\end{equation}
A typical situation when \eqref{coef} is satisfied is when the principal part  $A$ is uniformly elliptic and Lipschitz continuous in both  $x$ and  $t$.  Our main result which is the parabolic counterpart of the strong unique continuation result in \cite{Ru} can now  be stated as follows.

\begin{thrm} \label{main} 
Let $u \in L^{\infty}( B_2 \times [0,2))$ be a solution to \eqref{sub} in $B_2 \times [0, 2)  $  where $V$ satisfies
\begin{equation}\label{vthm}
||V||_{L^{\infty}(B_2 \times [0,2))} \leq M
\end{equation}
 and the coefficient matrix  $(a_{ij})$ satisfies the assumptions in \eqref{coef}. 
 
 Now if $u$ vanishes to infinite order in space in the sense of Definition \ref{van} below,  then we have that   $u(x,0)=0$ for all $x \in B_2$. 
\end{thrm} 

Similar to \cite{Ru}, our proof of Theorem \ref{main} is based on  the following   new Carleman  estimate for sublinear parabolic operators  which in turn  is based on a somewhat  delicate  adaptation of the techniques in  \cite{EF} to our sublinear situation.  As the reader will see, the proof of the following estimate is made possible by  combination of several non-trivial geometric facts which  thanks  to the specific structure of the sublinearity,  beautifully combine.  Moreover, unlike the elliptic case, the proof of Theorem \ref{main} following the Carleman estimate  is somewhat more involved because the ensuing inequalities are in the Gaussian space. 
\begin{thrm} \label{carlemanthm}
For a given $\alpha \geq 1$,  with  $\gamma=\frac{\alpha}{\delta^2}$,   let $u \in C_{0}^{\infty}( B_2 \times (0, \frac{1}{2\gamma}))$ be a solution to
\begin{equation}\label{e0}
\operatorname{div}(A\nabla u) +\partial_t u + f((x,t), u)= g, \text{$1< q <2$},
\end{equation}
where  the coefficient matrix $A=(a_{ij})$ satisfies \eqref{coef}. 
Define  $\sigma$ as in Lemma \ref{lma4} below  corresponding to $\theta$  as in \eqref{theta} and $\gamma$ as above. Also let $G$ be as in  Lemma \ref{lma5}. Then  there are numbers $\delta_0, N_0$ and $\tilde C$ depending on $\lambda, C_0, \beta$ in \eqref{coef} as well as the parameters in \eqref{a1} such  that for
$\alpha \geq \tilde C$ and $\delta \leq \delta_0,$ the
following inequality holds,
\begin{align}\label{Car1}
&\alpha \int_{\R^{n+1}_+} \sigma^{-\alpha} \frac{\theta(\gamma t)}{t} |u|^2 G dX+ \int_{\R^{n+1}_+} \sigma^{1-\alpha} \frac{\theta(\gamma t)}{t} |\nabla u|^2 G dX
\\
& + O(\alpha) \int_{\R^{n+1}_+} \sigma^{-\alpha} F(X, u) G dX
\notag
 \\
 &\leq  N_0 \int_{\R^{n+1}_+} \sigma^{1-\alpha} |g|^2 G dX
+  e^{N_0 \alpha} \gamma^{\alpha +N_0}\int_{\R^{n+1}_+}(u^2+t|\nabla u|^2 + F(X, u))  dX.
\notag
\end{align}
\end{thrm} 
In closing, we would like to mention that it remains to be seen whether one can also establish  a backward uniqueness result  for sublinear equations of the type  \eqref{sub} under global growth assumptions on the solution similar to that in \cite{Po}, \cite{Ch},  \cite{ESS} and \cite{WZ}.   It also seems to be a challenging open problem as to whether the regularity assumptions on the principal part  $A$ in Theorem \ref{main} can be further relaxed as in \cite{KT}.  We would like to address such questions in a future study.

\medskip

The paper is organized as follows.  In Section  \ref{2}, we introduce some basic notations and gather some known results that are relevant to our work. In Section  \ref{3}, we finally prove our main results. 

 \section{Notations and Preliminaries}\label{2}
In this section we introduce some basic notations and also  collect some background results from \cite{EF}  which will be used throughout our work.   
Given $r>0$ we  denote by $B_r(x_0)$ the Euclidean ball centered at  $x_0 \in \R^n$ and  when $x_0=0$, we denote it simply by $B_r$. From now on, a  generic  point $(x,t)$ in $\R^n \times [0, \infty)$  denoted by $X$.  Also, unless and otherwise  specified, $\nabla U, \operatorname{div}\ U$ will  refer to $\nabla_x U, \operatorname{div}_x U$ respectively.   The region  $\Rn \times (0, \infty)$  in space-time will be denoted by $\R^{n+1}_+$.  The notation $A \lesssim B$
would be mean  $A \leq CB$ for some universal $C$.  Also, for a  given a function  $\sigma(t)$, $\dot{\sigma} (t)$ would refer to its derivative. 

\medskip

We now state the relevant results from \cite{EF}. The first lemma  is a parabolic  Rellich type identity which  corresponds to Lemma $1$ in \cite{EF} and  similar to that in  \cite{EF} and \cite{EFV}, constitutes   the key  ingredient in the  proof of  our sublinear  Carleman estimate.  

\begin{lemma} \label{lma1}
Let $\sigma=\sigma(t)$ be a non-decreasing function satisfying $\sigma(0)=0, ~\alpha \in \R,$ and $H$ and $G$ denote two functions in $\R_{+}^{n+1},~G$ non-negative. Then, the following identity holds for all $u \in C_0^{\infty}(\R^{n+1}_+),$
\begin{align}
&2\int_{\R^{n+1}_+} \frac{\sigma^{1-\alpha}}{\dot{\sigma}}\left(\partial_t u- <A\nabla \log G , \nabla u> +\frac{1}{2} Hu-\frac{\alpha \dot{\sigma}}{2 \sigma}u\right)^2 G dX
\\
& + \int_{\R^{n+1}_+} \frac{\sigma^{1-\alpha}}{\dot{\sigma}} \partial_t[\log \frac{\sigma}{\dot{\sigma} t}] \, <A\nabla u, \nabla u> G dX
\notag
\\
&= 2\int_{\R^{n+1}_+} \frac{\sigma^{1-\alpha}}{\dot{\sigma}}(\operatorname{div}(A \nabla u) +\partial_t u)\left(\partial_t u-<A \nabla \log G ,\nabla u>+\frac{1}{2} Hu-\frac{\alpha \dot{\sigma}}{2 \sigma}u\right) G dX
\notag
\\
&+ \int_{\R^{n+1}_+} \frac{\sigma^{1-\alpha}}{\dot{\sigma}} u<A\nabla u, \nabla H> G dX-\frac{1}{2} \int_{\R^{n+1}_+} \frac{\sigma^{1-\alpha}}{\dot{\sigma}} u^2 M dX
\notag \\
&+\frac{\alpha}{2} \int_{\R^{n+1}_+} \sigma^{-\alpha} u^2(\partial_t G-\operatorname{div}(A\nabla G)-HG) dX
\notag
 \\
&-\int_{\R^{n+1}_+} \frac{\sigma^{1-\alpha}}{\dot{\sigma}} < A\nabla u, \nabla u> (\partial_t G- \operatorname{div} (A\nabla G)-HG)   dX-2\int_{\R^{n+1}_+} \frac{\sigma^{1-\alpha}}{\dot{\sigma}} D_G \nabla u \cdot \nabla u  dX
\notag
\end{align}

where
$$M=\partial_t[\log \frac{\sigma}{\dot{\sigma} t}] HG+\partial_t HG+H(\partial_t G-\operatorname{div}(A \nabla G)-HG)-<A\nabla G,\nabla H>;$$
$D_G$ is the $n \times n$ symmetric matrix defined as 
\begin{eqnarray*}
D_G^{ij}=\frac{a_{ij}}{2t} G+a_{il} \partial_{kl} Ga_{kj}-\frac{a_{ik}\partial_kG a_{jl} \partial_l G}{G}+\frac{1}{2} \partial_ka_{il} \partial_l G a_{kj}+\frac{1}{2} \partial_k a_{jl} \partial_l G a_{ki}-\frac{1}{2} a_{kl} \partial_l G \partial_k a_{ij}+\frac{1}{2} \partial_t a_{ij} G.
\end{eqnarray*}
\end{lemma}

We also need the following identity ( see Lemma 2 in \cite{EF}).

\begin{lemma} \label{lma2}
Assume that $\sigma$ and $G$ are as  in Lemma \ref{lma1}. Then, the following identity
holds for  $u \in C_0^{\infty}(\R^{n+1}_+)$ and $\alpha \in \R$,
\begin{eqnarray*}
(\alpha-1) \int_{\R^{n+1}_+} \sigma^{-\alpha} \partial_t[\log \frac{\sigma}{\dot{\sigma} t}] \, u^2 G dX
&=& 2\int_{\R^{n+1}_+} \frac{\sigma^{1-\alpha}}{\dot{\sigma}} \, \partial_t[\log \frac{\sigma}{\dot{\sigma} t}] \, [u( \operatorname{div}(A \nabla u)+\partial_t u)+|\nabla u|^2] G dX\\
&+& 2\int_{\R^{n+1}_+} \frac{\sigma^{1-\alpha}}{\dot{\sigma}} \partial_t[\log \frac{\sigma}{\dot{\sigma} t}] u^2(\partial_t G-\operatorname{div}(A \nabla G)) dX \\
&-&2\int_{\R^{n+1}_+} \sigma^{1-\alpha} \partial_t \left(\frac{1}{\dot{\sigma}}\partial_t[\log \frac{\sigma}{\dot{\sigma} t}]\right) u^2 G \, dX.
\end{eqnarray*}

\end{lemma}

As in \cite{EF}, $\sigma$ as in \eqref{Car1} is chosen to be a solution to an appropriate ordinary differential equation  which is  dictated by the identity above.  To this end, we have the following Lemma which is Lemma 4 in \cite{EF}. 

\begin{lemma} \label{lma4}
Assume that $\theta:(0,1) \to \R_+$ satisfies
$$0\leq \theta \leq N,~|t\dot{\theta}(t)| \leq N\theta(t) \mbox{ and } \int_0^1(1+\log \frac{1}{t}) \frac{\theta(t)}{t} \, dt \leq N$$
for some constant $N.$ Then the solution to the ordinary differential equation 
$$\frac{d}{dt} \log (\frac{\sigma}{t\dot{\sigma}})= \frac{\theta(\gamma t)}{t},~\sigma(0)=0,~\dot{\sigma}(0)=1,$$
where $\gamma >0,$ has the following properties when $0\leq \gamma t\leq 1$:
\begin{enumerate}
\item $t e^{-N} \leq \sigma(t) \leq t,$
\item $e^{-N} \leq \dot{\sigma}(t)\leq 1,$
\item $|\partial_t[\sigma \log \frac{\sigma}{\dot{\sigma} t}]|+|\partial_t[\sigma \log \frac{\sigma}{\dot{\sigma} }]|\leq 3N$,
\item $\left|\sigma \partial_t \left(\frac{1}{\dot{\sigma}}\partial_t[\log \frac{\sigma}{\dot{\sigma} t}]\right)\right| \leq 3N e^{N} \frac{\theta(\gamma t)}{t}.$
\end{enumerate}

\end{lemma}
Now corresponding to $\beta$ as in  \eqref{coef}, the function $\theta$ is chosen  as follows
\begin{equation}\label{theta}
 \theta(t)=t^{\beta/2}\left(\log \frac{1}{t}\right)^{1+\beta/2}.
 \end{equation}
 It is easily seen that $\theta$ satisfies the conditions in Lemma \ref{lma4}.   From now, let  $0 < \delta < 1$ denote a small number to be chosen later, and $\alpha$ and $\beta$ be 
two numbers satisfying $\alpha \geq  1$ and $0 < \beta \leq 1$. We also need the following weighted inequalities  in the Gaussian space ( see Lemma 5 in \cite{EF}). 
\begin{lemma} \label{lma5}
Let $G(X)=t^{-n/2} e^{-|x|^2/4t}$ and $\sigma$ denote the function defined in Lemma \ref{lma4} corresponding to  $\gamma=\frac{\alpha}{\delta^2}$ and $\theta$ as in \eqref{theta}.

Then, there is a constant $N$ depending on $\beta$ and $n$ such that the following inequalities hold for all functions $u \in C_{0}^{\infty}(\R^n \times [0, 1/2\gamma))$,

\begin{eqnarray*}
&&\int_{\R^{n+1}_+} \sigma^{-\alpha} u^2\left(\frac{|x|^{\beta}}{t}+\frac{|x|^{2+\beta}}{\alpha t^2}+t^{\beta/2 -1}\right) G dX\\
&&\leq N e^{N \alpha} \gamma^{\alpha +N}\int_{\R^{n+1}_+} u^2 dX+N\delta^{\beta} \int_{\R^{n+1}_+} \sigma^{-\alpha} \frac{\theta(\gamma t)}{t} u^2 G dX;\\
&&\int_{\R^{n+1}_+} \sigma^{1-\alpha} |\nabla u|^2\left(\frac{|x|^{\beta}}{t}+\frac{|x|^{2+\beta}}{\alpha t^2}+ \frac{|x|^{1+\beta}}{t \delta}+t^{\beta/2 -1}\right) G dX \\
&&\leq N e^{N \alpha} \gamma^{\alpha +N}\int_{\R^{n+1}_+}t|\nabla u|^2 dX+N\delta^{\beta} \int_{\R^{n+1}_+} \sigma^{1-\alpha} \frac{\theta(\gamma t)}{t} |\nabla u|^2 G dX.
\end{eqnarray*}

\end{lemma}

 In closing, we    define  the  relevant notion of vanishing to infinite order in space. 

\begin{dfn}\label{van}
We say that a function $u$  defined in a region $\Om$ in space time  vanishes to infinite order in space  at $(x_0, t_0) \in \Om$ if given $k>0$, there exists $C_k>0$ such
\begin{equation} \label{ist1.1}
|u(x,t_0)|\leq C_k|x-x_0|^k
\end{equation}
for all $(x, t_0) \in \Om$. 
\end{dfn}

\section{Proof of the main results Theorem \ref{carlemanthm} and  Theorem \ref{main}}\label{3}

\subsection*{Proof of Theorem \ref{carlemanthm}}
\begin{proof}
By rotation of coordinates, without loss of generality we may
assume that $A(0,0)=\mathbb{I}$. Then as in \cite{EF}, we let  $r(x) = |x|$ and 
\[
H = \frac{r^2(1- < A\nabla r, \nabla r>)}{4 t^2}.
\]
By a standard calculation we have 
$$\partial_t G-\operatorname{div}( A \nabla G)=\left(\frac{r^2(1-< A \nabla r, \nabla r>)}{4 t^2}+ \frac{( < A \nabla r, \nabla r>-1)}{2 t} +\frac{r \operatorname{div}( A \nabla r)-(n-1)}{2 t}\right) G$$
and 
\begin{eqnarray} \label{estF}
|H|\lesssim \text{min} \left(\frac{(|x|^2+ t)^{1+\beta/2}}{t^2}, \frac{|x|^2}{t^2} \right), ~\left|\partial_t G-\operatorname{div}( A \nabla G)-HG\right| \lesssim \frac{(|x|^2+ t)^{\beta/2}}{t} G.
\end{eqnarray}
Note that \eqref{estF} in particular implies that
\begin{equation}\label{est10}
\left|\partial_t G -\operatorname{div}( A \nabla G)\right| \lesssim  \left(\frac{(|x|^2+ t)^{1+\beta/2}}{t^2} + \frac{(|x|^2+ t)^{\beta/2}}{t}  \right) G.
\end{equation}

Now  as in the statement of  Theorem \ref{carlemanthm}, for a given $\alpha$,   we have
\[
\gamma=\alpha/\delta^2.
\]
Let   $\theta, \sigma$ be also  as in Theorem \ref{carlemanthm}. Then by  using    the identity in Lemma \ref{lma2}  and  the equation \eqref{e0} satisfied by $u$,  the estimates in  \eqref{est10} and  the bounds for $\sigma,\dot{\sigma}$  in Lemma \ref{lma4}, we get the following estimate
\begin{align}\label{c80}
&(\alpha-1) \int_{\R^{n+1}_+} \sigma^{-\alpha} \frac{\theta(\gamma t)}{t} u^2 G dX
\\
&  \leq  2\int_{\R^{n+1}_+} \frac{\sigma^{1-\alpha}}{\dot{\sigma}}  \, \frac{\theta(\gamma t)}{t} \, u(-f(X, u)+  g) G dX
\notag
\\
& + C [\int_{\R^{n+1}_+} \frac{\sigma^{1-\alpha}}{\dot{\sigma}} \frac{\theta(\gamma t)}{t} \left( |\nabla u|^2 +u^2  \left(\frac{(|x|^2+ t)^{1+\beta/2}}{t^2} + \frac{(|x|^2+ t)^{\beta/2}}{t}  \right)  \right) G dX
\notag
\\
&+ \int_{\R^{n+1}_+} \sigma^{-\alpha}\frac{\theta(\gamma t)}{t} u^2 G dX].
\notag
\end{align}
Note that  in \eqref{c80} above, we also used the differential equation satisfied by $\sigma$ as in Lemma \ref{lma4}.  Now  since $\frac{\theta(\gamma t)}{t} \lesssim \sigma^{-1}$, therefore  by  applying   the weighted inequality in Lemma \ref{lma5} to the term
\[
\int_{\R^{n+1}_+} \frac{\sigma^{1-\alpha}}{\dot{\sigma}} \frac{\theta(\gamma t)}{t} u^2  \left(\frac{(|x|^2+ t)^{1+\beta/2}}{t^2} + \frac{(|x|^2+ t)^{\beta/2}}{t}  \right) G dX
\]

 we deduce  that  the  following holds
\begin{align}\label{c1}
&(\alpha-1) \int_{\R^{n+1}_+} \sigma^{-\alpha} \frac{\theta(\gamma t)}{t} u^2 G dX
\\
&\leq C\big[ \int_{\R^{n+1}_+} \sigma^{1-\alpha} \, \frac{\theta(\gamma t)}{t} \, |\nabla u|^2 G dX  +  e^{N \alpha} \gamma^{\alpha +N}\int_{\R^{n+1}_+} u^2  dX +(\alpha \delta^{\beta}+1) \int_{\R^{n+1}_+} \sigma^{-\alpha} \frac{\theta(\gamma t)}{t} u^2 G dX  \big] 
  \notag
\\
& + 2\int_{\R^{n+1}_+} \frac{\sigma^{1-\alpha}}{\dot{\sigma}}  \, \frac{\theta(\gamma t)}{t} \, u(-f(X, u)+  g) G dX\notag
\end{align}
for some universal $C, N$ depending also  on  the bounds in \eqref{coef}.    Now  observe that if $\delta$ is small enough and $\alpha$ is taken large enough, then the following integral in \eqref{c1} 

\[
C( \alpha \delta^{\beta} +1) \int_{\R^{n+1}_+} \sigma^{-\alpha} \frac{\theta(\gamma t)}{t} u^2 G dX
\]
can be absorbed in the left hand side.  Moreover since
\[
uf(X, u) \geq 0,
\] 
the following integral on the right hand side of \eqref{c1} 
\[
2\int_{\R^{n+1}_+} \frac{\sigma^{1-\alpha}}{\dot{\sigma}} \, \frac{\theta(\gamma t)}{t} u(-f(X, u)) G
\]
is non-positive and hence the inequality in \eqref{c1} remains valid  without this term. 
 Then    by applying   Cauchy Schwartz inequality  to   
 \[
2 \int_{\R^{n+1}_+} \frac{\sigma^{1-\alpha}}{\dot{\sigma}} \, \frac{\theta(\gamma t)}{t}\  ug\ G dX 
 \]
we deduce from \eqref{c1}   that the following estimate holds, 
\begin{align} \label{mainest}
&\alpha \int_{\R^{n+1}_+} \sigma^{-\alpha} \frac{\theta(\gamma t)}{t} u^2 G dX
\\
&\leq C_1\big[ \int_{\R^{n+1}_+} \sigma^{1-\alpha} \, \frac{\theta(\gamma t)}{t} \, |\nabla u|^2 G dX+\int_{\R^{n+1}_+} \sigma^{1-\alpha} \, | g|^2 G dX
\notag
\\
 & + e^{N \alpha} \gamma^{\alpha +N}\int_{\R^{n+1}_+} u^2  dX  \big]
 \notag
\end{align}
for some $C_1$ universal.
Now in order to incorporate the gradient term on the left hand side  in the Carleman estimate, we make use of the identity in Lemma \ref{lma1}. 
For that,   we first  note that  using
\[
\frac{\delta_{ij}}{2t} G + \partial_{ij} G - \frac{\partial_i G \partial_j G}{G} =0 
\]
and the bounds on the derivatives of $A$ as in  \eqref{coef} that the following estimate holds, 
\begin{eqnarray} \label{estD}
|D_G \nabla u\cdot \nabla u| \lesssim \frac{(|x|^2 +t)^{\beta/2}}{t} \, |\nabla u|^2 G.
\end{eqnarray}
This  corresponds to the estimate  (3.3) in \cite{EF}. 
 Next from Lemma \ref{lma4}, the bounds on the coefficients as in   \eqref{coef}   and \eqref{estF},  we have 
\begin{eqnarray} \label{Mest}
|\sigma \nabla H| \lesssim \frac{(|x|^2+t)^{1+\beta/2}}{t} \mbox{ and } |\sigma M| \lesssim \frac{(|x|^2+t)^{1+\beta/2}}{t^2} G.
\end{eqnarray}
Then by applying the  identity as  in Lemma \ref{lma1} and by using  the equation  \eqref{e0} satisfied by $u$  we obtain

\begin{align}\label{r0}
&2\int_{\R^{n+1}_+} \frac{\sigma^{1-\alpha}}{\dot{\sigma}}\left(\partial_t u- < A\nabla \log G , \nabla u> +\frac{1}{2} Hu-\frac{\alpha \dot{\sigma}}{2 \sigma}u\right)^2 G dX
\\
& + \int_{\R^{n+1}_+} \frac{\sigma^{1-\alpha}}{\dot{\sigma}} \partial_t[\log \frac{\sigma}{\dot{\sigma} t}] \, < A\nabla u, \nabla u> G dX
\notag
\\
&= 2\int_{\R^{n+1}_+} \frac{\sigma^{1-\alpha}}{\dot{\sigma}}(-f(X, u) +  g)\left(\partial_t u-< A \nabla \log G ,\nabla u>+\frac{1}{2} Hu-\frac{\alpha \dot{\sigma}}{2 \sigma}u\right) G dX
\notag
\\
&+ \int_{\R^{n+1}_+} \frac{\sigma^{1-\alpha}}{\dot{\sigma}} u< A\nabla u, \nabla H> G dX-\frac{1}{2} \int_{\R^{n+1}_+} \frac{\sigma^{1-\alpha}}{\dot{\sigma}} u^2 M dX
\notag \\
&+\frac{\alpha}{2} \int_{\R^{n+1}_+} \sigma^{-\alpha} u^2(\partial_t G-\operatorname{div}( A\nabla G)-HG) dX
\notag
 \\
&-\int_{\R^{n+1}_+} \frac{\sigma^{1-\alpha}}{\dot{\sigma}} < A\nabla u, \nabla u> (\partial_t G- \operatorname{div} ( A\nabla G)-HG)   dX-2\int_{\R^{n+1}_+} \frac{\sigma^{1-\alpha}}{\dot{\sigma}} D_G \nabla u \cdot \nabla u  dX.
\notag
\end{align}

Now, using  Cauchy-Schwarz inequality, the  following integral  on the right hand side
\[
2 \int_{\R^{n+1}_+} \frac{\sigma^{1-\alpha}}{\dot{\sigma}} g \left(\partial_t u-< A \nabla \log G ,\nabla u>+\frac{1}{2} Hu-\frac{\alpha \dot{\sigma}}{2 \sigma}u\right) G dX
\]
can be estimated as 
\begin{align}\label{r1}
& 2 \int_{\R^{n+1}_+} \frac{\sigma^{1-\alpha}}{\dot{\sigma}} g \left(\partial_t u-< A \nabla \log G ,\nabla u>+\frac{1}{2} Hu-\frac{\alpha \dot{\sigma}}{2 \sigma}u\right) G dX
\\
& \leq \frac{1}{4} \int_{\R^{n+1}_+} \frac{\sigma^{1-\alpha}}{\dot{\sigma}} \left(\partial_t u-< A \nabla \log G ,\nabla u>+\frac{1}{2} Hu-\frac{\alpha \dot{\sigma}}{2 \sigma}u\right)^2  G dX
\notag
\\
& + 8 \int_{\R^{n+1}_+} \frac{\sigma^{1-\alpha}}{\dot{\sigma}} g^2 G dX
\notag
\end{align}
and then the first integral in  the right hand side of  \eqref{r1} can be absorbed into the first term in the left hand side of \eqref{r0}. 
 Consequently  it follows  from \eqref{r0} and \eqref{r1}  and by  using  the bounds \eqref{estF}, \eqref{estD}, \eqref{Mest} as well as the inequalities in  Lemma \ref{lma5} that the following estimate holds 
\begin{align} \label{imp2est}
& \int_{\R^{n+1}_+} \sigma^{1-\alpha} \, \frac{\theta(\gamma t)}{t} \, |\nabla u|^2 G dX  
\\
&\leq  \int_{\R^{n+1}_+} \frac{\sigma^{1-\alpha}}{\dot{\sigma}}(-f(X, u) )\left(\partial_t u- <A \nabla \log G , \nabla u>+\frac{1}{2} Hu-\frac{\alpha \dot{\sigma}}{2 \sigma}u\right) G dX  \notag\\
&+ C[ \int_{\R^{n+1}_+} \sigma^{1-\alpha} \, |g|^2 G dX 
+ \alpha \delta^{\beta}\int_{\R^{n+1}_+} \sigma^{-\alpha} \frac{\theta(\gamma t)}{t} |u|^2 G dX\notag \\ &
\delta^{\beta} \int_{\R^{n+1}_+} \sigma^{1-\alpha} \, \, \frac{\theta(\gamma t)}{t} |\nabla u|^2 G dX
+ \int_{\R^{n+1}_+} \sigma^{-\alpha}|u| |\nabla u| \frac{(|x|^2+t)^{(1+\beta)/2}}{t} G dX \nonumber
\\&+  e^{N \alpha} \gamma^{\alpha +N}\int_{\R^{n+1}_+}(u^2+t|\nabla u|^2)  dX].
\notag
\end{align}
Finally the term 
\[
\int_{\R^{n+1}_+} \sigma^{-\alpha}|u| |\nabla u| \frac{(|x|^2+t)^{(1+\beta)/2}}{t} G dX
\]
 is handled using   Cauchy-Schwarz inequality in the following way
\begin{eqnarray} \label{imp1est}
 \int_{\R^{n+1}_+} \sigma^{-\alpha}|u| |\nabla u| \frac{(|x|^2+t)^{(1+\beta)/2}}{t} G dX &&\lesssim  \int_{\R^{n+1}_+} \sigma^{-\alpha}|u|^2  \frac{(|x|^2+t)^{1+\beta/2}}{t^2} G dX \nonumber\\
 &+& \int_{\R^{n+1}_+} \sigma^{1-\alpha} |\nabla u|^2 \frac{(|x|^2+t)^{\beta/2}}{t} G dX.
\end{eqnarray}
Now the  terms on the right hand side of \eqref{imp1est} are  again estimated using the inequalities in Lemma \ref{lma5}  and consequently we deduce from \eqref{imp2est} that the following holds
\begin{align} \label{r4}
& \int_{\R^{n+1}_+} \sigma^{1-\alpha} \, \frac{\theta(\gamma t)}{t} \, |\nabla u|^2 G dX  
\\
&\leq  \int_{\R^{n+1}_+} \frac{\sigma^{1-\alpha}}{\dot{\sigma}}(-f(X, u) )\left(\partial_t u- <A \nabla \log G , \nabla u>+\frac{1}{2} Hu-\frac{\alpha \dot{\sigma}}{2 \sigma}u\right) G dX  \notag\\
&+ C \bigg[ \int_{\R^{n+1}_+} \sigma^{1-\alpha} \, |g|^2 G dX 
+ \alpha \delta^{\beta}\int_{\R^{n+1}_+} \sigma^{-\alpha} \frac{\theta(\gamma t)}{t} |u|^2 G dX\notag \\ &
\delta^{\beta} \int_{\R^{n+1}_+} \sigma^{1-\alpha} \, \, \frac{\theta(\gamma t)}{t} |\nabla u|^2 G dX
+  e^{N \alpha} \gamma^{\alpha +N}\int_{\R^{n+1}_+}(u^2+t|\nabla u|^2)  dX \bigg].
\notag
\end{align}
Now by  using the fact that $A(0,0)=\mathbb{I}$ and the  bounds on the derivatives  $A$ as in \eqref{coef}, we observe that  
\begin{equation}\label{r5}
A \nabla \log G=-(\frac{x}{2t}+\frac{ O(|x|^2 +t)^{(\beta+1)/2}}{2t}).
\end{equation}
 We  note that the first term in the  right hand side of  \eqref{r4} can be equivalently written as
\begin{eqnarray} \label{b0}
&&\int_{\R^{n+1}_+} \frac{\sigma^{1-\alpha}}{\dot{\sigma}}(-f(X, u))\left(\partial_t u-<A\nabla \log G , \nabla u>+\frac{1}{2} Hu-\frac{\alpha \dot{\sigma}}{2 \sigma}u\right) G dX\\
&=&\int_{\R^{n+1}_+} \frac{\sigma^{1-\alpha}}{\dot{\sigma}}(-f(X,u))\left(\frac{Zu}{2t}+\frac{1}{2} Hu-\frac{\alpha \dot{\sigma}}{2 \sigma}u\right) G dX, \nonumber
\end{eqnarray}
where 
\[
Z=  2t( -<A \nabla \log G, \nabla> +  \partial_t).
\]
Now because of \eqref{r5} as well as \eqref{coef}, it follows that
\begin{equation}\label{r7}
\begin{cases}
Z= <x, \nabla> + 2t \partial_t + O((|x|^2 + t)^{(\beta+1)/2}) \nabla_x
\\
\operatorname{div}_X Z= n + 2+  O(|x|^{\beta}).
\end{cases}
\end{equation}

Now we look at each individual term in the right hand side of \eqref{b0}. First we observe  that from the following identity
\[
ZF(X,u)= f(X,u) Zu + <Z, \nabla_X F>
\]
  the  first term on the right hand side of \eqref{b0} can be rewritten as
\begin{align} \label{r8}
&-\int_{\R^{n+1}_+} \frac{\sigma^{1-\alpha}}{\dot{\sigma}}f(X,u) \frac{Z u}{2t} G dX= -  \int_{\R^{n+1}_+}  \frac{\sigma^{1-\alpha}}{\dot{\sigma}} \frac{1}{2t}( ZF(X,u)G - <Z, \nabla_X F>G) dX.
\end{align}
Now from the bounds in \eqref{a1} we see  that the second term in the right hand side  of \eqref{r8} can be upper bounded by 
\[
C_2 \int_{\R^{n+1}_+} \sigma^{1-\alpha}  \frac{1}{2t} (|x|+ |t|) F(X, u) G dX.
\]
Then again by using the first inequality in Lemma \ref{lma5} with $u$ replaced by $\sqrt{F}$, we  can assert that this term can by bounded from above in the following way
\begin{align}\label{r10}
&\int_{\R^{n+1}_+} \sigma^{1-\alpha} \frac{1}{2t}  (|x|+ |t|) F(X, u) G dX \leq  \int_{\R^{n+1}_+} \sigma^{1-\alpha} F(X,u) G dX 
\\
&+ C \left[  e^{N\alpha} \gamma^{\alpha+N}  \int_{\R^{n+1}_{+}}   t F(X,u) dX  + \delta^{\beta} \int_{\R^{n+1}_{+}} \sigma^{1-\alpha} \frac{\theta(\gamma t)}{t}  F(X,u)G dX \right].
\notag
\end{align}
In \eqref{r10}, we also used the fact that since $u(\cdot, t)$ is supported in $B_2$, therefore $|x| \lesssim |x|^{\beta}$ in the support of $u$. 

Now  by applying integration by parts to  the first integral in the right hand side of \eqref{r8}  we obtain
\begin{align}\label{f0}
& -\int_{\R^{n+1}_{+}}  \frac{\sigma^{1-\alpha}}{\dot{\sigma}} \frac{1}{2t} ZF(X,u)G dX
\\
& = \int_{\R^{n+1}_{+}}   \left[ \frac{\sigma^{1-\alpha}}{\dot{\sigma}} \frac{1}{2t}   F(X,u) ( (\operatorname{div}_X Z) G  + ZG)  +  F(X, u) G  (t \partial_t )[ \frac{\sigma^{1-\alpha}}{ t\dot{\sigma}}]  \right] dX.
\notag
\end{align}
At this point, we   note that  since $G$ is parabolic homogeneous of  degree $-n$, therefore we have that
\[
<x,  \nabla G> + 2t \partial_t G = - nG.
\]
Then by using this fact, it follows from the expression of $Z$ as in \eqref{r7} that the following holds, 
\begin{equation}\label{f2}
ZG=-(n+\frac{O((|x|^2 + t)^{1+\beta/2})}{t})G.
\end{equation}
We  also have that
\begin{equation}\label{f4}
 (t\partial_t)[\frac{\sigma^{1-\alpha}}{t\dot{\sigma}}]=(1-\alpha) \sigma^{-\alpha}-\frac{\sigma^{1-\alpha}}{t\dot{\sigma}}-\frac{\sigma^{1-\alpha}\ddot{\sigma}}{(\dot{\sigma})^2}. 
 \end{equation}
Now  note that since
$\frac{d}{dt}\text{log}(\frac{\sigma}{t \dot{\sigma}})= \frac{\dot{\sigma}}{\sigma}-\frac{1}{t}-\frac{\ddot{\sigma}}{\dot{\sigma}}=\frac{\theta(\gamma t)}{t}$, therefore
\begin{equation}
\frac{\ddot{\sigma}}{\dot{\sigma}}\lesssim \frac{1}{t}.
\end{equation}
This implies  that 
\begin{equation}\label{f5}
 (t\partial_t)[\frac{\sigma^{1-\alpha}}{t\dot{\sigma}}] =  -\alpha \sigma^{-\alpha} + O(1) \sigma^{-\alpha}.
 \end{equation}
 
  Therefore by using  \eqref{r7},  \eqref{f2} and \eqref{f5} in \eqref{f0}  it follows
\begin{align}\label{f7}
&-\int_{\R^{n+1}_{+}}  \frac{\sigma^{1-\alpha}}{\dot{\sigma}} \frac{1}{2t} ZF(X,u)G dX
\\
&= \int_{\R^{n+1}_+} \left[ \frac{\sigma^{1-\alpha}}{\dot{\sigma}} \left( \frac{1}{t} +\frac{O((|x|^2 + t)^{1+\beta/2})}{t} + \frac{O(|x|^{\beta})}{t} \right)G F(X,u) +(  - \alpha + O(1)) \sigma^{-\alpha} F(X, u) G \right]dX.
\notag
\end{align}
Now the first term on the right hand side of \eqref{f7} is again estimated by using the inequalities in Lemma \ref{lma5} ( with $\sqrt{F}$ instead of $u$) as follows
\begin{align}\label{f10}
& \int_{\R^{n+1}_+} \frac{\sigma^{1-\alpha}}{\dot{\sigma}} \left( \frac{1}{t} +\frac{O((|x|^2 + t)^{1+\beta/2})}{t} + \frac{O(|x|^{\beta})}{t} \right)G F(X,u) dX
\\
& \leq C \bigg[ e^{N\alpha} \gamma^{\alpha+N} \int_{\R^{n+1}_+}   F(X, u) dX +  \delta^{\beta} \int_{\R^{n+1}_+} \sigma^{1-\alpha} \frac{\theta(\gamma t)}{t} F(X, u) G dX 
\notag
\\
&+  \int_{\R^{n+1}_+} \sigma^{-\alpha} F(X, u) G dX \bigg].
\notag
\end{align}
Over here we note that in \eqref{f10} above,  we used the fact that in the support of $u(\cdot, t)$ which is contained in $B_2$,  we have 
\[
\frac{O((|x|^2+ t)^{1+\beta/2}}{t} \lesssim \frac{O(|x|^{\beta})}{t}  + O(1).
\]
We also used the bounds for $\sigma, \dot{\sigma}$ as in Lemma \ref{lma4}.   Now  again by using the bounds for $\sigma, \dot{\sigma}$ as in Lemma \ref{lma4}, the bounds for  $H$ as in \eqref{estF}  and the fact that
\[
0\leq u f(X,u) \leq q F(X, u)
\]
which is contained in the structural assumptions  as in   \eqref{a1} ( In fact this is precisely the place where we use the specific structure of the sublinearity),  we obtain the following estimate for the last two terms in the right hand side of \eqref{b0}, 
 \begin{align}\label{b8}
 &| \int_{\R^{n+1}_+} \frac{\sigma^{1-\alpha}}{\dot{\sigma}}(-f(X,u))\left(\frac{1}{2} Hu-\frac{\alpha \dot{\sigma}}{2 \sigma}u\right) G dX  |
 \\
 & \leq  C \int_{\R^{n+1}_+}  \sigma^{1-\alpha} \frac{ (|x|^2+t)^{1+\beta/2}}{t^2} F(X,u) GdX + \frac{\alpha q}{2} \int_{\R^{n+1}_+} \sigma^{-\alpha} F(X, u) GdX.
 \notag
 \end{align}
 Now again by using the inequalities in Lemma \ref{lma5}, the first term in the right hand side of \eqref{b8} can be estimated in the following way
 \begin{align}\label{b10}
 & \int_{\R^{n+1}_+}  \sigma^{1-\alpha} \frac{ (|x|^2+t)^{1+\beta/2}}{t^2} F(X,u) GdX 
 \\
 & \leq  N e^{N\alpha} \gamma^{\alpha+N} \int_{\R^{n+1}_+} F(X, u) dX + N \delta^{\beta} \alpha \int_{\R^{n+1}_{+}} \sigma^{1-\alpha} \frac{\theta(\gamma t)}{t} F(X, u) G dX.
 \notag
 \end{align} 
 
 At this point, by using the estimates \eqref{r10}, \eqref{f7}, \eqref{f10}, \eqref{b8} and \eqref{b10} in \eqref{r4} we obtain
 \begin{align} \label{r11}
& \int_{\R^{n+1}_+} \sigma^{1-\alpha} \, \frac{\theta(\gamma t)}{t} \, |\nabla u|^2 G dX  
\\
&\leq  (-\alpha+\frac{q}{2} \alpha + O(1)) \int_{\R^{n+1}_+} \sigma^{-\alpha} F(X,u) G dX \notag\\
&+ C\bigg[ \int_{\R^{n+1}_+} \sigma^{1-\alpha} \, |g|^2 G dX 
+ (\alpha \delta^{\beta} +1) \int_{\R^{n+1}_+} \sigma^{-\alpha} \frac{\theta(\gamma t)}{t} u^2 G + \sigma^{-\alpha} F(X,u) G dX\notag \\ &
+\delta^{\beta} \int_{\R^{n+1}_+} \sigma^{1-\alpha} \, \, \frac{\theta(\gamma t)}{t} |\nabla u|^2 G dX
+  e^{N \alpha} \gamma^{\alpha +N}\int_{\R^{n+1}_+}(u^2+t|\nabla u|^2 + F(X,u))  dX \bigg].
\notag
\end{align}
Over here we note that in  order to get to \eqref{r11} above, we also used the fact that $\frac{\theta(\gamma t)}{t} \lesssim c \sigma^{-1}$. Now the inequality  \eqref{mainest} can be equivalently written as

\begin{align} \label{m1}
&\frac{\alpha}{C_1} \int_{\R^{n+1}_+} \sigma^{-\alpha} \frac{\theta(\gamma t)}{t} u^2 G dX
\\
&\leq \big[ \int_{\R^{n+1}_+} \sigma^{1-\alpha} \, \frac{\theta(\gamma t)}{t} \, |\nabla u|^2 G dX+\int_{\R^{n+1}_+} \sigma^{1-\alpha} \, | g|^2 G dX
\notag
\\
 & e^{N \alpha} \gamma^{\alpha +N}\int_{\R^{n+1}_+} u^2  dX  \big].
 \notag
\end{align}
Now since $q < 2$, we can now choose $\alpha$ sufficiently large such that 
\begin{equation}\label{t1}
 (-\alpha+\frac{q}{2} \alpha + O(1)) \leq \frac{1}{4} (q-2)\alpha .
 \end{equation}
 
 Moreover, we  also choose $\delta$  small enough such that
 \begin{equation}\label{t2}
 C\delta^{\beta} < \frac{1}{16}, \  C(\alpha \delta^{\beta}+1) \leq \text{min}( \frac{\alpha}{ 4 C_1}, \frac{1}{8} (2-q) \alpha)
 \end{equation}
 where $C$ and $C_1$ are  the constants as in \eqref{r10} and  \eqref{m1} respectively.  Now by adding the inequalities \eqref{r10} $+$ $\frac{1}{2}$ \eqref{m1}, the desired estimate in Theorem \ref{carlemanthm} follows by also taking into account \eqref{t1} and \eqref{t2}.

\end{proof}
Before proceeding further, we make the following remark.

\begin{rmrk}
It remains to be seen whether the Carleman estimate \eqref{Car1} holds if  we  instead assume that $A$ satisfies 
\begin{equation}\label{st}
|A(x,t) - A(y,s)| \leq C (|x-y|^2 + |t-s|)^{1/2}
\end{equation}
as  in Theorem 1 i) in  \cite{EF}  ( which constitutes an alternate set of conditions under which the strong unique continuation result in \cite{EF} is valid)  or more generally if $A$ can be  allowed to have   $1/3$ H\"older regularity  in time as in \cite{KT}. However   the proof of the Carleman estimate in \cite{EF} for  principal part $A$ with regularity assumptions as in \eqref{st} crucially  relies on weighted Calderon Zygmund estimate of the following type
\begin{align}
& \frac{1}{\delta^2} \int \sigma^{2-\alpha} |\nabla^2u|^2 G dX  \lesssim \int \sigma^{1-\alpha} |\operatorname{div}(A(x,0) \nabla u)|^2 G dX + e^{N \alpha} \gamma^{\alpha+N} \int (u^2 + t |\nabla u|^2)G dX
\end{align}
( see for instance (3.9) in \cite{EF}).
It remains to be seen whether in  our sublinear situation, one can get similar estimates with $|\operatorname{div}(A(x,0) \nabla u)|^2$ replaced by $(\operatorname{div}( A(x,0) \nabla u) + f(X, u))^2$. This appears to be a challenging interesting issue  to which we would like to come back in a future study. 

\end{rmrk}

\subsection*{Proof of Theorem \ref{main}}

 We now proceed with the proof of our main unique continuation result Theorem \ref{main}.

\begin{proof}
The proof is divided into two steps. 

\medskip

\emph{Step 1}: We first assume that $u$  vanishes to infinite order in space and time at $(0,0)$, i.e. for every $k \in \mathbb{N}$, there exists $C_k$ such that 
\begin{equation}\label{v}
|u(x,t)| \leq C_k ( |x|^2+t)^{k/2}.
\end{equation}
Also  by taking a smaller neighborhood if necessary and then by iteratively spreading the zero set, without loss of generality we may assume that $|u| \leq 1$.  Now note that from our regularity assumption on $A$ as in \eqref{coef}, it follows from the  Calderon-Zygmund estimates as in \cite{Li} that  given any $p< \infty$,  

\begin{equation}\label{reg}
 \nabla^2 u, u_t \in L^{p} (B_r \times [0, 2))
\end{equation}
for all $r< 2$. Moreover, we also have from the Schauder theory   that $\nabla u$ is in $H^{\beta}_{loc}(B_2 \times [0, 2))$. 
   Now for a given $\ve>0$, let $u_\ve=u \phi_{\epsilon}(t) \, \psi(x),$ where $\psi \in C^{\infty}_0(\R^n)$  satisfies $\psi=1$ for $|x|\leq 1,~\psi=0$ for $|x| \geq 3/2$ and  $\phi_\ve\in C^{\infty}_0(\R)$ is a smooth cutoff such that
\begin{equation}
\begin{cases}
\phi_\ve \equiv 1,\ \text{when $\epsilon \leq t \leq \frac{1}{4\gamma}$}
\\
\phi_\ve \equiv 0, \ \text{when $t\leq \frac{\epsilon}{2}$ or $t \geq \frac{1}{2\gamma}.$}
\end{cases}
\end{equation}

Then we have that
\begin{eqnarray*}
&&\sum_{i,j}^n \partial_i(a_{ij}(X)\partial_ju_{\epsilon}) +\partial_tu_{\epsilon}+ f(X, u_\ve)=\tilde{g}_\ve
\end{eqnarray*}
where \begin{eqnarray*}
&&\tilde{g}_{\epsilon}=f(X, u_\ve) -f(X, u) \phi_{\epsilon}(t) \psi(x) -Vu \phi_\ve(t) \psi(x)\\
&&+2\phi_{\epsilon}(t) \sum_{i,j}^n a_{ij}(X)\partial_ju \partial_i \psi
+u\phi_{\epsilon}(t) \sum_{i,j}^n \partial_i(a_{ij}(X)\partial_j\psi)+u \phi_{\epsilon}^{'}(t) \psi(x).
\end{eqnarray*}

Now  given an integer $k\geq  \tilde C$ (where $\tilde C$ is as in Theorem \ref{carlemanthm})    we apply the Carleman estimate  in Theorem \ref{carlemanthm}
with $\alpha=2k$ to $u_{\epsilon}$ ( note that the validity of the Carleman estimate  \eqref{Car1} for $u_\ve$ can be justified using an approximation  with smooth functions and also by using   \eqref{reg}). Consequently we have

\begin{align}\label{t}
&\alpha \int_{\R^{n+1}_+} \sigma^{-\alpha} \frac{\theta(\gamma t)}{t} |u_{\epsilon}|^2 G dX+ \int_{\R^{n+1}_+} \sigma^{1-\alpha} \frac{\theta(\gamma t)}{t} |\nabla u_{\epsilon}|^2 G dX 
\\
&+ O(\alpha) \int_{\R^{n+1}_+} \sigma^{-\alpha} F(X, u_\ve) G dX \notag\\ &\leq N_0\int_{\R^{n+1}_+} \sigma^{1-\alpha} |\tilde{g}_{\epsilon}|^2 G dX
+  e^{N_0 \alpha} \gamma^{\alpha +N_0}\int_{\R^{n+1}_+}(u_{\epsilon}^2+t|\nabla u_{\epsilon}|^2+ F(X, u_\ve))  dX.
\notag
\end{align}

Now using $uf(X,u) \leq q F(X, u)$ we note   that the following Caccioppoli  type energy estimate hold
\begin{equation}\label{d1}
\int_{B_r \times (a, a+ r^2)} |\nabla u|^2 dX \leq \frac{C}{r^2}  \int_{B_{2r} \times (a, a+ 2r^2)} (  u^2 + F(X,u))  dX.
\end{equation}
Now since 
\[
F(X, u) \leq c |u|^{p}
\]
for some $p\in (1,2)$, it follows from \eqref{v}, the gradient estimate above  and parabolic regularity estimates as in \cite{Li} that $\nabla u$ also vanishes to infinite order in space and time in the sense of \eqref{v}.

Now by splitting  the integrals  over $\R^n \times (0, \infty)$   into dyadic time-like  regions of the type $\R^n \times \{1/2^{k} \leq t \leq 1/2^{k-1} \}$   and by using vanishing to infinite order property of $u$  and $\nabla u$, we can assert that for any $\alpha>0$, 
\[
\int_{B_{3/2} \times (0, 1)} \sigma^{-\alpha} ( u^2 + |\nabla u|^2)G < \infty.
\]
Moreover using the following bound
\[
|\phi_{\ve}'(t)| \leq \frac{C}{\ve}, \  \text{when $\frac{\ve}{2} < t < \ve$}
\]
and \eqref{v}, we note that  as $\ve \to 0$, 
\[
\int_{\frac{\ve}{2} < t< \ve}  u^2 \psi^2 \phi_{\ve}'(t)^2 dX \to 0.
\]

 Therefore we  can let $\ve \to 0$ in \eqref{t} and consequently we  obtain for $u_0= u \phi_0 \psi$ where $\phi_0$ is the pointwise limit of $\phi_\ve$ that the following inequality holds,

\begin{align}\label{4}
& \alpha \int_{\R^{n+1}_+} \sigma^{-\alpha} \frac{\theta(\gamma t)}{t} |u_{0}|^2 G dX+ \int_{\R^{n+1}_+} \sigma^{1-\alpha} \frac{\theta(\gamma t)}{t} |\nabla u_{0}|^2 G dX + O(\alpha) \int_{\R^{n+1}_+} \sigma^{-\alpha} F(X, u_0) G dX
\\
&\leq N_0\int_{\R^{n+1}_+} \sigma^{1-\alpha} |\tilde{g}_{0}|^2 G dX + e^{N_0 \alpha} \gamma^{\alpha +N_0}\int_{\R^{n+1}_+}(u_{0}^2+t|\nabla u_{0}|^2+|u_{0}|^q)  dX
\notag
\\
& \leq N_0\int_{\R^{n+1}_+} \sigma^{1-\alpha}  [(f(X, u) \phi_0 \psi - f(X, u_0))^2  + |u|^2 \phi_0^2 \psi^2] G  dX\notag\\
&+  2N_0\int_{B_2\setminus B_1 \times (0,1/2\gamma)} \sigma^{1-\alpha} |\phi_{0}(t) \nabla u|^2 G dX+N_0\int_{B_2\setminus B_1 \times (0,1/2\gamma)} \sigma^{1-\alpha} |\phi_{0}(t) u|^2 G dX
\notag
\\
&+N_0\int_{\R^{n+1}_+} \sigma^{1-\alpha}|u \phi_{0}^{'}(t) \psi|^2 G dX+ e^{N_0 \alpha} \gamma^{\alpha +N_0}\int_{\R^{n+1}_+}(u_{0}^2+t|\nabla u_{0}|^2+|u_{0}|^p)  dX
\notag
\end{align}
where $N_0$ additionally depends on the $L^{\infty}$ norm of $V$.

Now we estimate   each individual term in the right hand side of the above expression. We first  note that from the expression of $\theta$ as in  Lemma \ref{lma5} it follows that  $\gamma \lesssim \frac{\theta(\gamma t)}{t}$. Consequently using $\gamma= \frac{\alpha}{\delta^2}$, we  have
\begin{equation}\label{gamma}
 \delta^{-2} \lesssim \frac{\theta(\gamma t)}{t}
\end{equation}
 when $0< 2\gamma t <1.$
Therefore we have that the following   terms on the right hand side of \eqref{4} 
\[
2N_0\int_{B_2\setminus B_1 \times (0,1/2\gamma)} \sigma^{1-\alpha} |\phi_{0}(t) \nabla u|^2 G dX+N_0\int_{B_2\setminus B_1 \times (0,1/2\gamma)} \sigma^{1-\alpha} |\phi_{0}(t) u|^2 G dX  + N_0 \int_{\R^{n+1}_+} \sigma^{1-\alpha} u^2 \phi_0^2 \psi^2 G dX
\]
 can be estimated in the following way 

\begin{align}\label{c5}
& 2N_0\int_{B_2\setminus B_1 \times (0,1/2\gamma)} \sigma^{1-\alpha} |\phi_{0}(t) \nabla u|^2 G dX+N_0\int_{B_2\setminus B_1 \times (0,1/2\gamma)} \sigma^{1-\alpha} |\phi_{0}(t) u|^2 G dX 
\\
& +  N_0 \int_{\R^{n+1}_+} \sigma^{1-\alpha} u^2 \phi_0^2 \psi^2 G dX
\notag
\\
& \leq C \delta ^{2} \left( \int_{B_2  \times (0,1/2\gamma)} \sigma^{1-\alpha} \frac{\theta(\gamma t)}{t} |\nabla u|^2 G dX +
 \int_{B_2  \times (0,1/2\gamma)}  \sigma^{1-\alpha} \frac{\theta(\gamma t)}{t}|u|^2 G dX \right)
\notag
\end{align}
where  $C$ is independent of $\delta$. 
Now by choosing $\delta$ sufficiently small, we note that these terms can then be absorbed in the left hand side of  \eqref{4}. We consequently fix such a $\delta$.  We then consider the following term  in the right hand side of \eqref{4} 
\[
 \int_{\R^{n+1}_+} \sigma^{1-\alpha}  (f(X, u) \phi_0 \psi - f(X, u_0))^2 G dX. 
 \]
 Note that this term is non-zero  in $B_2 \setminus B_1 \times (0, 1/2\gamma) \cup  B_2 \times (1/4\gamma, 1/2\gamma)$.
 Therefore this term can be estimated from above in the following way using $e^{-N} t \leq \sigma \leq t$, 
 \begin{align}\label{c0}
 & \int_{\R^{n+1}_+} \sigma^{1-\alpha}  \left(f(X, u) \phi_0 \psi - f(X, u_0)\right)^2 G dX 
 \\
 & \leq C e^{2Nk} \left( \int_{B_{3/2} \setminus B_1 \times (0, 1/2\gamma)}  t^{-2k+1} |u|^{2p-2} G dX + \int_{B_{3/2} \times (1/4\gamma, 1/2\gamma)} t^{-2k+1} |u|^{2p-2} G dX \right)
 \notag
 \end{align}
 
 where  in \eqref{c0} above, we used that $|f(\cdot, s)| \leq \kappa |s|^{p-1}$ for some $p \in (1, 2)$.  Now since $\alpha= 2k$ and $\gamma= \frac{\alpha}{\delta^2}$, it follows that if $|x|\geq 1$ or $t \geq 1/4\gamma$, there exists $N$ 
 depending also on $\delta$ such that
 \[ 
 t^{-2k} G\leq e^{2Nk} k^{2k}. 
 \]
 Also from Stirling's formula, we have 
 \[
 k^{k} \lesssim e^{Nk} k!.
 \]
 
 Therefore for a new $N$, we obtain
 
 \begin{align}\label{c2}
 & \int_{\R^{n+1}_+} \sigma^{1-\alpha}  (f(X, u) \phi_0 \psi - f(X, u_0))^2 G dX 
 \\
 & \leq C e^{2Nk} (k!)^2 \left( \int_{B_{3/2} \setminus B_1 \times (0, 1/2\gamma)}   |u|^{2p-2}  dX + \int_{B_{3/2} \times (1/4\gamma, 1/2\gamma)}  |u|^{2p-2} dX \right).
 \notag
 \end{align} 
 Likewise the term 
 \[
 N_0\int_{B_2\setminus B_1 \times (0,1/2\gamma)} \sigma^{1-\alpha} |\phi_{0}(t) u|^2 G dX
 \]
 can be bounded from above  in the following way
 \begin{align}\label{c4}
 & N_0\int_{B_2\setminus B_1 \times (0,1/2\gamma)} \sigma^{1-\alpha} |\phi_{0}(t) u|^2 G dX  \leq
 C e^{2Nk} (k!)^2  \int_{B_{3/2} \setminus B_1 \times (0, 1/2\gamma)}   u^2  dX.
 \end{align}
 
 Then we observe that the term 
\[
N_0\int_{\R^{n+1}_+} \sigma^{1-\alpha}|u \phi_{0}^{'}(t) \psi|^2 G dX
\]
is estimated  from above as follows

\begin{align}\label{c10}
& N_0\int_{\R^{n+1}_+} \sigma^{1-\alpha}|u \phi_{0}^{'}(t)|^2 G dX\leq N\int_{B_2\times (1/4\gamma, 1/2\gamma)} \sigma^{1-\alpha}|u|^2  G dX\lesssim e^{2Nk} (k!)^2  \int_{B_2\times (1/4\gamma, 1/2\gamma)}  u^2  dX
\end{align}
where in \eqref{c10} above, we again  made use of Stirling formula.  Finally the following  term  in the right hand side of \eqref{4}
\[
e^{N_0 \alpha} \gamma^{\alpha +N_0}\int_{\R^{n+1}_+}(u_{0}^2+t|\nabla u_{0}|^2+|u_{0}|^p)  dX
\]
is handled as follows
\begin{align}
& e^{N_0 \alpha} \gamma^{\alpha +N_0}\int_{\R^{n+1}_+}(u_{0}^2+t|\nabla u_{0}|^2+|u_{0}|^p)  dX \leq e^{2Nk} (k!)^2 \int_{B_{3/2} \times (0, 1/2\gamma)} (u^2 + t |\nabla u|^2 + |u|^p) dX.
\end{align}
Now by using the energy estimate as in  \eqref{d1}, we can assert that the following inequality holds,
\begin{align}\label{c7}
& e^{N_0 \alpha} \gamma^{\alpha +N_0}\int_{\R^{n+1}_+}(u_{0}^2+t|\nabla u_{0}|^2+|u_{0}|^p)  dX \leq e^{2Nk} (k!)^2 \int_{B_{2} \times (0, 2)}  |u|^p dX
\end{align}
for some $N$. Note that in order to get to \eqref{c7}, we used the  Stirling formula  and also the fact that since $u$ is bounded, therefore $u^2 \lesssim |u|^p$.  Therefore, by combining \eqref{4}, \eqref{c5},  \eqref{c2}, \eqref{c4}, \eqref{c10}  and \eqref{c7} we finally obtain for a new $N$ and $k \geq \tilde C$ that the following holds, 

\begin{align}\label{c89}
& \int_{B_1 \times (0, 4/\gamma)} t^{-2k} u^2 GdX \leq (\frac{1}{4} N)^{2+2k}( k!)^2  \int_{B_{2} \times (0,2)} |u|^{2(p-1)} dX.
\end{align}
Note that in \eqref{c89} we also used the boundedness of $u$ and the fact that since $p \in (1,2)$, therefore 
\[
2(p-1)= \text{min} ( 2(p-1), p, 2).
\]

Now   by writing 
\[
\int_{B_1 \times (0,1)} t^{-2k} u^2 G dX = \int_{B_1 \times (0,4/\gamma)} t^{-2k} u^2 G dX  + \int_{B_1 \times (4/\gamma, 1)} t^{-2k} u^2 G dX
\]
and by estimating 
\[
\int_{B_1 \times (4/\gamma, 1)} t^{-2k} u^2 G dX
\]
using
\[
t^{-2k} G \lesssim e^{2Nk}  k^{2k}, \ \text{since $t\geq 4/\gamma$}
\]

we consequently obtain  using Stirling formula and \eqref{c89} that for a new $N$,  the following estimate holds for $k \geq \tilde C$, 
\begin{align}\label{c8}
&\left( \int_{B_1 \times (0, 1)} t^{-2k} u^2 GdX \right)^{1/2} \leq (\frac{1}{4} N)^{1+k}( k!)  ||u||^{(p-1)}_{L^{\infty}(B_{2} \times (0, 2))}. 
\end{align}
Now by multiplying the inequality  by $2^k/N^k k!$ and summing over $k \geq \tilde C$, we obtain
\begin{align}\label{d2}
& \sum_{k\geq \tilde C} \left(\int_{B_1 \times (0,1)} \frac{2^k}{N^k t^k k!} u^2 G dX \right)^{1/2} \leq N ||u||^{p-1}_{L^{\infty}(B_{2} \times (0, 2))}. 
\end{align}
Now we note that there exists $K_0$ depending on $\tilde C$, such that  for $a \geq K_0$, 
\begin{equation}\label{d4}
\sum_{k \geq \tilde C} \frac{a^k}{k!} \geq \frac{1}{2} e^a.
\end{equation}
Consequently we have from  \eqref{d2}, \eqref{d4} and triangle inequality that the following holds, 
\begin{align}\label{d5}
&  \left(\int_{B_1 \times (0,\frac{2}{NK_0})} e^{2/Nt} u^2 G dX \right)^{1/2} \lesssim ||u||^{p-1}_{L^{\infty}(B_{2} \times (0, 2))}. 
\end{align}
Now  using 
\[
\frac{1}{Nt} \geq \frac{|x|^2}{8t} - \frac{|x-y|^2}{8t}
\]
when $|y| \leq 8/N$ and $|x| \leq 1/2$,  we obtain from \eqref{d5} and the explicit expression of $G(x, t)$   that
\begin{align}\label{d7}
&  \left(\int_{B_1 \times (0,s)} t^{-n/2} u^2 e^{-\frac{|x-y|^2}{4t}} dX \right)^{1/2} \lesssim e^{-1/Ns} ||u||^{p-1}_{L^{\infty}(B_{2} \times (0, 2))} 
\end{align}
for  $s \in  (0, \frac{1}{NK_0})$ and $y$ such that $|y| \leq 8/N$. Therefore it follows from \eqref{d7} that for all  such $s, y$ we have 

\begin{equation}\label{d8}  
\left(\int_{B_{\sqrt{s}}(y) \times (s, 2s)} s^{-n/2} u^2  dX \right)^{1/2} \lesssim  e^{-1/Ns} ||u||^{p-1}_{L^{\infty}(B_{2} \times (0, 2))}. 
\end{equation}

Now note that since $u$ solves \eqref{sub}, therefore by treating $f(X,u)$ as a scalar term and by applying the standard Moser subsolution estimate as  in \cite{AS}( see also  Theorem 6.29 in \cite{Li}), we obtain
\begin{equation}\label{d10}
|u(y,s)| \lesssim \frac{1}{s^{n/2 +1}} \int_{s}^{2s} \int_{B_{\sqrt{s}}(y)}  |u| dX +  s^{1- n/2q_0} \left(\int_{s}^{2s} \int_{B_{\sqrt{s}}(y)} u^{(p-1) q_0} dX \right)^{1/q_0}, \ \text{$q_0> n/2+1$}.
\end{equation}
Here we also used that $f(\cdot, s) \lesssim |s|^{p-1}$ as in \eqref{a0}. 
We  additionally choose $q_0$ large enough such that $(p-1)q_0  \geq 2$.
Now by Cauchy-Schwartz and \eqref{d8}, the first integral  in the right hand side of \eqref{d10} is upper bounded by $ e^{-1/Ns} ||u||^{(p-1)/2}_{L^{\infty}(B_{2} \times (0, 2))}$ and for the  second integral, since $(p-1)q_0  \geq 2$ and $u$ is bounded, therefore the second integral on the right hand side  of \eqref{d10} can be estimated as follows
\begin{equation}\label{d11}
s^{1- n/2q_0} \left(\int_{s}^{2s} \int_{B_{\sqrt{s}}(y)} u^{(p-1) q_0} dX \right)^{1/q_0} \leq C \left(\int_{s}^{2s} \int_{B_{\sqrt{s}}(y)} u^{2} dX \right)^{1/q_0}
\end{equation}
where $C$ depends on the  $L^{\infty}$ norm of $u$ which again  because of \eqref{d8} can be upper bounded  by $e^{-1/Ns} ||u||^{2(p-1)/q_0}_{L^{\infty}(B_{2} \times (0, 2))}$ for a  different $N$. Therefore finally we obtain that for some universal $N$ that $u$ satisfies the following estimate
\begin{equation}\label{g1}
|u(y, s) | \leq C e^{-1/Ns}, \ \text{when $|y| \leq  8/N$, $s < \frac{1}{NK_0}$}. 
\end{equation}
This implies that $u(\cdot, 0) \equiv 0$ in $B_{8/N}$  and  the estimate \eqref{g1} in particular implies that  $u$ vanishes to infinite order in space and time at every $(y, 0)$ for $|y| \leq 8/N$. At this point, by a standard argument  we can spread the zero set  and conclude that $u(\cdot, 0) \equiv 0$. 

\medskip

\emph{Step 2}: We now show that if $u$ vanishes to infinite order at $(0,0)$ in the space variable in the sense of  Definition \ref{van}, then $u$ also vanishes to infinite order in both space and time  in the sense of \eqref{v}.  For linear parabolic equations, this follows from a result of Alessandrini and Vessella in \cite{AV}. We note that  proof in \cite{AV} uses the local asymptotics of solutions to parabolic equations vanishing to a certain order in space and time as derived in \cite{AV1}. However the proof of such a local asymptotic result in \cite{AV1} relies on   certain  scaling properties of a linear equation  in a crucial way and this is  not available in our  sublinear situation.  Therefore  we instead adapt an alternate approach due to Fernandez in \cite{F}. 

\medskip

We proceed as follows.
We note that it suffices to show \eqref{v} for $k$ large enough. Let $k \geq M$ where $M$ is large enough to be decided later.  We additionally assume that  $M \geq \tilde C$ where $\tilde C$ is as in Theorem \ref{carlemanthm}. Corresponding to this $k$,  as before let $\alpha=2k$ and $\gamma=\frac{\alpha}{\delta^2}$ where $\delta$ is small enough as required in \emph{Step} $1$. 
Now  for a fixed  $a$ such that  $0< a < \frac{1}{4\gamma}$ and   with  $\psi, \phi_0$  as in \emph{Step 1} corresponding to such a  $\gamma$,  by repeating the arguments in the proof of Theorem \ref{carlemanthm}  to $u_0= u \psi \phi_0$ in the region $\R^{n+1}_+$ with
$G(x, t+a)$ instead of $G(x, t)$  and $\sigma(t+a)$ instead of $\sigma(a)$, and  by keeping track of the additional positive boundary terms which occur when integrating by parts with respect to the time-variable and then by adding up such terms to the right hand side of  our previous estimate  \eqref{4}, we note that  after such a computation, the additional boundary integrals   on the right hand side ( i.e. at $\R^n \times \{0\}$)  are  bounded  from above by a multiple of  
\begin{equation}\label{h0}
\alpha \sigma(a)^{-\alpha} \int_{\Rn \times \{0\}} u_0^2 (x, 0) \left( 1+ \frac{|x|^2}{a} \right)  G(x, a) dx\  +  \sigma(a)^{-\alpha} \int_{\Rn \times \{0\}} F((x,0), u_0(x, 0)) G(x, a) dx.
\end{equation}
We note that the first integral   above is as in \cite{F} ( see Section $3$ in \cite{F}) whereas the second integral  is the one  that is incurred due to an integration by parts of an expression involving the sublinear term as in \eqref{f0}.  Then by using
\[
F(\cdot, s) \leq C s^{p}, \text{for some $p \in (1, 2)$}
\]
we see that the  expression in \eqref{h0} is upper bounded by
\[
C\alpha \sigma(a)^{-\alpha} \int_{\Rn \times \{0\}} |u_0|^p (x, 0) \left( 1+ \frac{|x|^2}{a}  \right)G(x, a)  dx.
\]
Now by repeating the arguments  as in \emph{Step 1}  upto \eqref{c8} we obtain the following estimate for some universal $N$  ( Over here, note that the inequality  \eqref{gamma} still holds since $a \leq 1/4\gamma$)

\begin{align}\label{e10}
& \int_{B_1 \times (0, 1)} (t+a)^{-2k} u^2 G(x, t+a)dX \leq  N^{2k}( k!)^2  ||u||^{2(p-1)}_{L^{\infty}(B_{2} \times (0, 2))}  +  N^{2k} a^{-2k}   \int_{B_{2} \times \{0\} } |u|^p(x, 0) G(x,a) dx.
\end{align}

Now since $a \leq \frac{1}{4\gamma}$ and $\gamma \sim k$, therefore, we have that $a \leq \frac{1}{ Ck}$. Now note that \eqref{e10} in particular implies the following estimate
\begin{align}\label{h1}
& \int_{B_{\sqrt{a}} \times (0, a)}  u^2 dX \leq  N^{2k}( k!)^2 a^{2k} ||u||^{2(p-1)}_{L^{\infty}(B_{2} \times (0, 2))}  +  N^{2k}    \int_{B_{2} \times \{0\} } |u|^p(x, 0) G(x,a) dx.
\end{align}

 Now given some $a \in (0, \frac{1}{Ck})$, using the fact that $u$ vanishes to infinite order in space, we can ensure that 
\begin{equation}\label{e11}
N^{2k}   \int_{B_{2} \times \{0\} } |u|^p(x, 0) G(x,a) dx \leq M_k a^k.
\end{equation}
Now again by repeating the arguments as in \eqref{d10}-\eqref{d11} which uses the Moser's subsolution estimate and also by using \eqref{e11} we can assert  that there exists universal constants  $N_1, N_2>0$ such that for all $k \geq N_1$, we have  
\begin{equation}\label{h2}
||u||_{L^{\infty}(B_{\sqrt{a}/2} \times [0, a/2))} \leq  M_k a^{k/N_2}, \ \text{ for some $M_k$ and where $a \leq \frac{1}{Ck}$}.
\end{equation}
Now it can be seen by a standard real analysis argument that  \eqref{h2} implies \eqref{v} and consequently    by \emph{Step 1},  we can again  conclude  that $u(\cdot, 0) \equiv 0$. This completes the proof of the Theorem.

\end{proof}

\end{document}